\documentclass[a4paper,10pt]{article}

\usepackage{microtype}
\usepackage{color}
\usepackage[utf8x]{inputenc}
\usepackage{amsfonts}
\usepackage{amsmath}
\usepackage{amssymb}
\usepackage{amsthm}
\usepackage{mathrsfs}
\usepackage{latexsym}
\usepackage{graphicx}
\usepackage{bm}
\usepackage{inputenc}
\usepackage{enumerate}
\usepackage{enumitem}
\usepackage{hyperref}

\swapnumbers

\newtheorem{lemma}{Lemma}[section]
\newtheorem{theorem}[lemma]{Theorem}
\newtheorem{corollary}[lemma]{Corollary}

\theoremstyle{definition}
\newtheorem{definition}[lemma]{Definition}

\theoremstyle{remark}
\newtheorem*{Proof}{Proof}
\newtheorem{remark}[lemma]{Remark}

\newtheoremstyle{citing}
  {3pt}
  {3pt}
  {\itshape}
  {}
  {\bfseries}
  {.}
  {.5em}
  {\thmnote{#3}}
\theoremstyle{citing}
\newtheorem*{citing}{}

\DeclareMathOperator{\with}{:}
\DeclareMathOperator{\without}{\sim}
\DeclareMathOperator{\restrict}{\llcorner}
\DeclareMathOperator{\Clos}{Clos}  
\DeclareMathOperator{\Tan}{Tan}

\DeclareMathOperator{\Lip}{Lip}     
  
\DeclareMathOperator{\dmn}{dmn}     
\DeclareMathOperator{\Nor}{Nor}
     
\DeclareMathOperator{\Der}{D}       

  \DeclareMathOperator{\ap}{ap}
\DeclareMathOperator{\diam}{diam}
\DeclareMathOperator{\dist}{dist} 
\DeclareMathOperator{\Dis}{Dis}

\DeclareMathOperator{\aff}{af{}f}

\begin{document}

\title{A geometric second-order-rectifiable stratification for
closed subsets of Euclidean space} \author{Ulrich Menne \and Mario Santilli}

\maketitle

\begin{abstract}

	Defining the $m$-th stratum of a closed subset of an $n$~dimensional
	Euclidean space to consist of those points, where it can be touched by
	a ball from at least $n-m$ linearly independent directions, we
	establish that the $m$-th stratum is second-order rectifiable of
	dimension~$m$ and a Borel set.  This was known for convex sets, but is
	new even for sets of positive reach.  
	The result is based on a new
	criterion for second-order rectifiability.
\end{abstract}

\paragraph{MSC-classes 2010.} 52A20 (Primary); 28A78, 49Q15 (Secondary)


\section{Introduction}

The main purpose of the present paper is to establish the following theorem;
our notation is based on \cite[pp.~669--676]{MR41:1976}, see the end of this
introduction.

\begin{citing} [Structural theorem on the singularities of closed sets
$\mathrm{(see~\ref{thm:stratification-borel})}$] \hypertarget{structural}
	Suppose $A$ is a closed subset of $\mathbf R^n$, for $a \in A$, $\Dis
	(A,a)$ is the set of $v \in \mathbf R^n$ satisfying $A \cap \{ x \with
	|x-(a+v)|<|v| \} = \varnothing$, $m$ is an integer, $0 \leq m \leq n$,
	and
	\begin{equation*}
		B = A \cap \{ a \with \textup{$\Dis(A,a)$ contains at least
		$n-m$~linearly independent vectors} \}.
	\end{equation*}

	Then, $B$ can be $\mathscr H^m$ almost covered by the union of a
	countable collection of $m$~dimensional, twice continuously
	differentiable submanifolds of~$\mathbf R^n$.
\end{citing}

In the terminology of \cite[p.~2]{arXiv:1701.07286v1} for $m \geq 1$, the
conclusion asserts that $B$ is countably $(\mathscr H^m, m)$~rectifiable of
class~$2$.  If $A$ is convex, then $B$ consists of the set of points, where
the dimension of the normal cone of~$A$ is at least $n-m$, see
\ref{remark:alberti}.  Hence, our theorem contains Alberti's structural
theorem on the singularities of convex sets, see~\cite[Theorem~3]{MR1384392}.
We also prove, that $B$~is a countably $m$~rectifiable Borel set, see
\ref{thm:stratification-borel}; in particular, if $m \geq 1$, then $B$ can be
covered (without exceptional set) by a countable family of images of
Lipschitzian functions from~$\mathbf R^m$ into~$\mathbf R^n$, and, if $m=0$,
then $B$ is countable.

Our approach rests on two pillars.  The first may be stated as follows.

\begin{citing} [Parametric criterion for second-order rectifiability
$\mathrm{(see~\ref{corollary:rectifiability_criterion})}$]
\hypertarget{criterion}
	Suppose $W$ is an $\mathscr L^n$ measurable subset of $\mathbf R^n$,
	$m$ is an integer, $1 \leq m \leq n$, $f : W \to \mathbf R^\nu$ is a
	locally Lipschitzian map, $Z = \mathbf R^\nu \cap \{ z \with \mathscr
	H^{n-m} ( f^{-1} [ \{ z \} ] ) > 0 \}$, and, for $\mathscr H^m$~almost
	all $z \in Z$, there exists an $m$~dimensional subspace~$U$
	of~$\mathbf R^\nu$ satisfying
	\begin{equation*}
		\limsup_{y \to x} |y-x|^{-2} \dist (f(y)-f(x), U) < \infty
		\quad \text{whenever $x \in f^{-1} [ \{ z \} ]$}.
	\end{equation*}

	Then, $Z$ can be $\mathscr H^m$ almost covered by the union of a
	countable collection of $m$~dimensional, twice continuously
	differentiable  submanifolds of~$\mathbf R^\nu$.
\end{citing}

Notice that $f^{-1} [ \{ z \} ]$ abbreviates $\{ x \with f(x) =z \}$, see
below.  The key to reduce this
\hyperlink{criterion}{criterion} to the
nonparametric case is the construction (in \ref{thm:para-coarea}) of a
countable collection $G$ of $m$~rectifiable subsets~$P$ of~$W$ with $\mathscr
H^m ( Z \without f [ \bigcup G ] ) = 0$ such that, for each $P \in G$, the
restriction $f|P$ is univalent and $(f|P)^{-1}$ is Lipschitzian.  The
nonparametric case was comprehensively studied in \cite{arXiv:1701.07286v1};
however, for the present purpose, also \cite{MR2472179} would be sufficient
(see \ref{remark:little_effort}).

The second pillar of the proof of the structural theorem is the next result
that we state here for the special case of a convex set $A$.  It concerns the
relation of the nearest point projection, $\boldsymbol \xi_A$, with the
tangent and normal cones of~$A$.

\begin{citing} [A geometric observation for convex sets
$\mathrm{(see~\ref{lemma:relative_interior}\
with~\ref{lemma:convex}\,\eqref{item:convex:reach}\,\eqref{item:convex:normal}, \ref{remark:dis_xi})}$] \hypertarget{observation}
	If $A$ is a nonempty closed convex subset of $\mathbf R^n$, $m$ is an
	integer, $1 \leq m < n$, $x \in \mathbf R^n \without A$, $a =
	\boldsymbol \xi_A (x)$, $\dim \Nor (A,a) \geq n-m$, $U$ is an $m$
	dimensional subspace of~$\mathbf R^n$, $U \subset \Tan (A,a)$, and
	$x-a$ belongs the relative interior of $\Nor (A,a)$, then
	\begin{equation*}
		\limsup_{y \to x} |y-x|^{-2} \dist ( \boldsymbol \xi_A
		(y)-a,U)< \infty.
	\end{equation*}
\end{citing}

This criterion and its generalisation to closed sets
in~\ref{lemma:relative_interior} owe much to Federer's treatment of sets of
positive reach (a concept that embraces convex sets and submanifolds of
class~$2$) in \cite{MR0110078}.  Since it is elementary, that the set $B$ in
the structural theorem is countably $m$~rectifiable, the parametric criterion
for second-order rectifiability then is readily applied with $f = \boldsymbol
\xi_A |W$ for suitable~$W$.

\paragraph{Connection to curvature measures} Instead of using second-order
rectifiability properties, curvature properties can also be studied via
general Steiner formulae.  This approach was taken, for sets of positive reach
and various more general classes of sets, by Federer in~\cite{MR0110078},
Stach{\'o} in~\cite{MR534512}, Z{\"a}hle in~\cite{MR849863}, Rataj and
Z{\"a}hle in~\cite{MR1846894}, and Hug, Last, and Weil in~\cite{MR2031455}; in
fact, \cite{MR534512} and \cite{MR2031455} treat arbitrary closed subsets of
Euclidean space.  Accordingly, the natural question (under investigation by
the second author) arises to characterise the relation of both notions of
curvature.

\paragraph{Connection to varifold theory}  The original motivation of the
first author for the present study was to create a deeper understanding of a
relation proven by Almgren in his area-mean-curvature characterisation of the
sphere in~\cite{MR855173}.  There, an equation relating the curvature measures
(similar to those of~\cite{MR849863}) of the convex hull of the support of a
certain varifold to the perpendicular part of the mean curvature of the
varifold is established in~\cite[\S\,6\,(2)]{MR855173}.  The results of
present paper shall serve as tools for further investigations of both authors
of the second-order rectifiability properties of classes of varifolds.

\paragraph{Acknowledgements} The first author thanks the participants of an
online reading seminar of~\cite{MR855173} for their early interest in these
developments.  The material of this paper originates from the PhD thesis of
the second author, supervised by the first author, submitted at the University
of Potsdam.  The paper was written while both authors worked at the Max Planck
Institute for Gravitational Physics (Albert Einstein Institute) and the
University of Potsdam.

\paragraph{Notation}

Our notation and terminology is that of \cite[pp.~669--676]{MR41:1976}, except
that, as in \cite[p.~8]{MR0370454}, we denote the image of $A$ under a
relation $r$ by
\begin{equation*}
	r [A] = \{ y \with \textup{$(x,y) \in r$ for some $x \in A$} \}.
\end{equation*}

\section{Coarea formula}

The purpose of the present section is to prove the
parametric criterion for second-order rectifiability
in~\ref{corollary:rectifiability_criterion}.  We begin by establishing a
theorem that allows to construct univalent parametrisations from a
Lipschitzian given one.

\begin{theorem} \label{thm:para-coarea}
	
	Suppose $W$ is an $\mathscr L^n$ measurable subset of $\mathbf R^n$,
	$m$ is an integer, $1 \leq m \leq n$, and $f: W \rightarrow
	\mathbf R^\nu $ is a locally Lipschitzian map.  Then, there exists a
	countable collection $G$ of compact subsets $P$ of $W$, such that
	$f|P$ is univalent and $(f|P)^{-1}$ is Lipschitzian, satisfying
	\begin{equation*}
		\textstyle \mathscr H^m \big ( \mathbf R^\nu \cap \{ z \with
		\mathscr H^{n-m} ( f^{-1} [\{z\}] ) > 0 \} \without \bigcup
		\{ f [P] \with P \in G \} \big ) = 0.
	\end{equation*}
	Moreover, each member of $G$ is contained in some $m$ dimensional
	affine plane.
\end{theorem}

\begin{Proof}
	We firstly consider \emph{the special case that $W$ is a compact
	subset of~$\mathbf R^n$}.  Choose $F : \mathbf R^n \to \mathbf R^\nu$
	with $F|W = f$ and $\Lip F = \Lip f < \infty$
	by Kirszbraun's theorem
	\cite[2.10.43]{MR41:1976}; in particular, $\Der F$ is a Borel function
	whose domain is a Borel set by \cite[3.1.2]{MR41:1976}.  Defining $Z_i
	= \mathbf R^\nu \cap \{ z \with \mathscr H^{n-m} ( f^{-1} [ \{ z \}
	] ) \geq 1/i \}$ whenever $i$ is a positive integer, we note that
	\begin{equation*}
		\textstyle \mathbf R^\nu \cap \{ z \with \mathscr H^{n-m} (
		f^{-1} [ \{ z \} ] ) > 0 \} = \bigcup_{i = 1}^\infty Z_i.
	\end{equation*}
	Moreover, the sets $Z_i$ are Borel sets by 
	\cite[2.10.26]{MR41:1976} and $(\mathscr H^m, m)$~rectifiable by
	\cite[3.2.31]{MR41:1976}.  We define, for every positive integer~$i$,
	the class $\Omega_i$ to consist of all families $G$ of compact subsets
	$P$ of $f^{-1} [ Z_i]$ such that
	\begin{equation*}
		\text{$f [ P ] \cap f [Q] \neq \varnothing$ if and only if $P
		= Q$}
	\end{equation*}
	whenever $P, Q \in G$, and such that
	\begin{gather*}
		\mathscr H^m (P) > 0, \quad \text{$f|P$ is univalent}, \quad
		\text{$(f|P)^{-1}$ is Lipschitzian}, \\
		\text{$P$ is contained in some $m$ dimensional affine
		subspace of $\mathbf R^n$}
	\end{gather*}
	whenever $P \in G$.  Clearly, each member of $\Omega_i$ is countable.
	Using Hausdorff's maximal principle (see \cite[p.~33]{MR0370454}), we
	choose maximal elements $G_i$ of $\Omega_i$.  The proof of the
	present case will be concluded by establishing
	\begin{equation*}
		\textstyle \mathscr H^m \big ( Z_i \without \bigcup \{ f [ P
		] \with P \in G_i \} \big ) = 0 \quad \text{for every positive
		integer $i$}.
	\end{equation*}

	For this purpose, fix such $i$ and define Borel sets $T = Z_i \without
	\bigcup \{ f [ P ] \with P \in G_i \}$ and $S = f^{-1} [T]$.
	If $T$ had positive $\mathscr H^m$ measure, then, noting
	\cite[2.10.35]{MR41:1976},
	\begin{equation*}
		\textstyle B = S \cap \big \{ w \with \| \bigwedge_m \Der
		F(w) \| > 0 \big \}
	\end{equation*}
	would be a Borel set and have positive $\mathscr L^n$ measure by the
	coarea formula \cite[3.2.22\,(3)]{MR41:1976} with $W$, $Z$, and $f$
	replaced by $S$, $T$, and $f|S$, since
	\begin{equation*}
		( \mathscr L^n \restrict S, n ) \ap \Der (f|S)(w) = \Der F(w)
		\quad \text{for $\mathscr L^n$ almost all $w \in S$}
	\end{equation*}
	by \cite[2.10.19\,(4)]{MR41:1976}.
	
	Consequently, identifying $\mathbf R^n \simeq \mathbf R^m \times
	\mathbf R^{n-m}$, there would exist
	a linear isometry $g :
	\mathbf R^n \to \mathbf R^n$ such that $\mathscr L^n (A) > 0$ with
	\begin{align*}
		A & = \textstyle B \cap \big \{ w \with \bigwedge_m ( \Der
		F(w) | g [ \mathbf R^m \times \{ 0 \} ] ) \neq 0 \big \} \\
		& = \textstyle g \big [ g^{-1} [B] \cap \{ x \with
		\bigwedge_m ( \Der (F \circ g)(x) | \mathbf R^m \times \{
		0 \} ) \neq 0 \} \big ]
	\end{align*}
	and, as $A$ would be a Borel set, $\eta \in \mathbf R^{n-m}$ so that
	$\mathscr L^m (R) > 0$ with
	\begin{align*}
		R = \mathbf R^m \cap \{ \xi \with ( \xi,\eta ) \in g^{-1}
		[A] \}
	\end{align*}
	by Fubini's theorem, see \cite[2.6.2\,(3)]{MR41:1976}.
	Since $R$ would be a Borel set, we could apply \cite[3.2.2]{MR41:1976}
	to the function $h : \mathbf R^m \to \mathbf R^\nu$ defined by $h(\xi)
	= (F \circ g) ( \xi, \eta )$ for $\xi \in \mathbf R^m$, and use the
	Borel regularity of $\mathscr H^m$ to construct a subset $P$ of $g [ R
	\times \{ \eta \} ]$ with $G_i \cup \{ P \} \in \Omega_i$, contrary to
	the maximality of $G_i$.

	To treat the \emph{general case}, we choose an increasing sequence of
	compact subsets $K_i$ of $\mathbf R^n$ with $\mathscr L^n \big ( W
	\without \bigcup_{i=1}^\infty K_i \big ) = 0$.
	Since, in conjunction with \cite[2.4.5]{MR41:1976},
	\cite[2.10.25]{MR41:1976} applied with $A$ replaced by $W \without
	\bigcup_{i=1}^\infty K_i$ implies that
	\begin{equation*}
		\textstyle \mathscr H^{n-m} \big ( f^{-1} [ \{ z \} ]
		\without \bigcup_{i=1}^\infty K_i \big ) = 0 \quad \text{for
		$\mathscr H^m$ almost all $z \in \mathbf R^\nu$}
	\end{equation*}
	and $\lim_{i \to \infty} \mathscr H^{n-m} (f^{-1} [ \{ z \} ] \cap
	K_i ) = \mathscr H^{n-m} ( f^{-1} [ \{ z \} ] )$ for such
	$z$, we readily infer the conclusion.
\end{Proof}

\begin{remark}
	The contradiction argument is inspired by \cite[3.2.21]{MR41:1976}.
\end{remark}

\begin{remark}
	For the nearest point projection onto a set of positive reach, the
	idea of exhaustion by means of images from lower dimensional parts of
	the domain of~$f$ is employed in \cite[4.15\,(3)]{MR0110078}.  The
	important additional feature of members $P$ in our collection $G$ is
	the Lipschitz continuity of $(f|P)^{-1}$.
\end{remark}

\begin{remark}
	One readily verifies that \ref{thm:para-coarea} also holds with $m=0$,
	but this will not be needed in the present paper.
\end{remark}

The \hyperlink{criterion}{parametric criterion for second-order
rectifiability} now reads as follows.  

\begin{corollary} \label{corollary:rectifiability_criterion}
	Under the hypotheses of~\ref{thm:para-coarea}, if
	\begin{equation*}
		Z = \mathbf R^\nu \cap \{ z \with \mathscr H^{n-m} ( f^{-1} [
		\{ z \} ] ) > 0 \},
	\end{equation*}
	and, for $\mathscr H^m$~almost all $z \in Z$, there exists an
	$m$~dimensional subspace~$U$ of~$\mathbf R^\nu$ satisfying
	\begin{equation*}
		\limsup_{y \to x} |y-x|^{-2} \dist ( f(y)-f(x), U ) < \infty
		\quad \text{whenever $x \in f^{-1} [ \{ z \} ]$},
	\end{equation*}
	then $Z$ can be $\mathscr H^m$~almost covered by the union of a
	countable collection of $m$~dimensional submanifolds of $\mathbf
	R^\nu$ of class~$2$.
\end{corollary}

\begin{Proof}
	Whenever $P \in G$, as $(f|P)^{-1}$ is Lipschitzian, we notice that,
	for $\mathscr H^m$ almost all $z \in Z \cap f[P]$, there exists an
	$m$~dimensional subspace $U$ of $\mathbf R^\nu$ such that
	\begin{equation*}
		\limsup_{f [P] \owns \zeta \to z} |\zeta-z|^{-2} \dist
		(\zeta-z, U) < \infty.
	\end{equation*}
	Therefore, the conclusion follows from
	\cite[5.3]{arXiv:1701.07286v1} and \cite[3.1.15]{MR41:1976}.
\end{Proof}

\begin{remark} \label{remark:little_effort}
	With little additional effort, the final argument could have been
	based on \cite[A.1]{MR2472179} instead of
	\cite[5.3]{arXiv:1701.07286v1} and \cite[3.1.15]{MR41:1976}.
\end{remark}

\begin{remark} \label{remark:coarea_tangent}
	In conjunction with the preceding corollary, the following observation
	will be useful.  \emph{If $B$ is a countably $(\mathscr H^m,m)$
	rectifiable subset of $\mathbf R^\nu$, then, for $\mathscr H^m$ almost
	all $b \in B$, there exists an $m$~dimensional subspace~$U$
	of~$\mathbf R^\nu$ such that $U \subset \Tan (B,b)$;} in fact,
	\cite[2.1.4, 3.1.21]{MR41:1976} reduce the problem to Borel sets $B$,
	in which case \cite[2.10.19\,(4), 3.2.17, 3.2.18]{MR41:1976} apply.
\end{remark}

\section{Convex sets}

In the present section, we mainly collect some basic properties of convex sets
and related definitions in \ref{def:dist}--\ref{thm:hausdorff_distance} for
convenient reference.  Additionally, we note an observation concerning convex
cones in \ref{def:cone}--\ref{corollary:cone_control}.

\begin{definition} \label{def:dist}
	Suppose $A \subset \mathbf R^n$ and $x \in \mathbf R^n$.  Then, the
	\emph{distance of~$x$ to~$A$} is denoted by $\dist (x,A) = \inf \{
	|x-a| \with a \in A \}$.
\end{definition}

\begin{remark} \label{remark:distance}
	If $A \neq \varnothing$, then $\dist ( \cdot, A )$ is real valued and
	$\Lip \dist (\cdot,A) \leq 1$.
\end{remark}

\begin{remark} \label{remark:distance_relation}
	If $R = ( \mathbf R^n \times A ) \cap \{ (x,a) \with |x-a| = \dist
	(x,A) \}$, then, using~\ref{remark:distance}, one verifies that $\{ a
	\with \textup{$(x,a) \in R$ for some $x \in B$} \}$ is bounded
	whenever $B$ is a bounded subset of $\mathbf R^n$.  Moreover, if $A$
	is closed, so is $R$.
\end{remark}

\begin{definition} [see \protect{\cite[4.1]{MR0110078}}] \label{def:xi}
	Suppose $A \subset \mathbf R^n$ and $U$ is the set of all $x \in
	\mathbf R^n$ such that there exists a unique $a \in A$ with $|x-a| =
	\dist (x,A)$.  Then, the \emph{nearest point projection onto~$A$} is
	the map $\boldsymbol \xi_A : U \to A$ characterised by the requirement
	$| x- \boldsymbol \xi_A (x)| = \dist (x,A)$ for $x \in U$.
\end{definition}

\begin{remark} \label{remark:xi}
	Using \ref{remark:distance_relation}, we obtain that the function
	$\boldsymbol \xi_A$ is continuous.  Moreover, if $A$ is closed, then
	$\dmn \boldsymbol \xi_A$ is a Borel set; in fact, one verifies, by
	means of \ref{remark:distance_relation}, that the function mapping $x
	\in \mathbf R^n$ onto $d(x) = \diam \{ a \with (x,a) \in R \} \in
	\overline{\mathbf R}$ is upper semicontinuous, and $\dmn \boldsymbol
	\xi_A = \{ x \with d(x)=0 \}$.
\end{remark}

\begin{definition} [see \protect{\cite[p.~xix]{MR3155183}}]
	If $A \subset \mathbf R^n$, then $\aff A$ denotes the \emph{affine
	hull} of $A$.
\end{definition}

\begin{definition} [see \protect{\cite[p.~7, p.~xx]{MR3155183}}]
	Suppose $C$ is a convex subset of $\mathbf R^n$.  Then, the
	\emph{dimension} of $C$, denoted by $\dim C$, is defined to be the
	dimension of $\aff C$, and the \emph{relative boundary}
	[\emph{interior}] of $C$ is defined to be the boundary [interior] of
	$C$ relative to $\aff C$.
\end{definition}

\begin{remark} \label{remark:relative_interior}
	If $V$ is the relative interior of $C$, then $V$ is convex, $\dim V =
	\dim C$, and
	\begin{equation*}
		c+t(v-c) \in V \quad \text{whenever $v \in V$, $c \in C$, and
		$0 < t \leq 1$};
	\end{equation*}
	in fact, reducing to the case $\aff C = \mathbf R^n$, this is
	\cite[1.1.9, 1.1.10, 1.1.13]{MR3155183}.
\end{remark}

\begin{lemma} \label{lemma:convex}
	Suppose $C$ is a nonempty closed convex subset of $\mathbf R^n$.

	Then, the following four statements hold.
	\begin{enumerate}
		\item \label{item:convex:reach} There holds $\dmn \boldsymbol
		\xi_C = \mathbf R^n$ and $\Lip \boldsymbol \xi_C \leq 1$.
		\item \label{item:convex:tangent} If $c \in C$, then $\Tan
		(C,c) = \mathbf R^n \cap \{ u \with \textup{$u \bullet v \leq
		0$ for $v \in \Nor (C,c)$} \}$ and
		\begin{equation*}
			C \subset \{ c + u \with u \in \Tan (C,c) \} \subset
			\aff C;
		\end{equation*}
		in particular, $\dim C = \dim \Tan (C,c)$.
		\item \label{item:convex:normal} If $c \in C$, then
		\begin{gather*}
			\Nor (C,c) = \{ v \with \boldsymbol \xi_C (c+v)=c \} =
			\mathbf R^n \cap \{ v \with \textup{$v \bullet (x-c)
			\leq 0$ for $x \in C$} \}.
		\end{gather*}
		\item \label{item:convex:boundary} If $B$ is the relative
		boundary of $C$, then
		\begin{gather*}
			B = C \cap \{ c \with \textup{$c+v \in \aff C$ for some
			$v \in \mathbf S^{n-1} \cap \Nor (C,c)$} \}.
		\end{gather*}
	\end{enumerate}
\end{lemma}

\begin{Proof}
	\eqref{item:convex:reach} is asserted in \cite[4.1.16]{MR41:1976}.  In
	view of~\eqref{item:convex:reach}, the first equation and the first
	inclusion in~\eqref{item:convex:tangent} are contained in
	\cite[4.8\,(12)]{MR0110078} and \cite[4.18]{MR0110078}, respectively;
	the remaining items of~\eqref{item:convex:tangent} then follow.  The
	first equation in~\eqref{item:convex:normal} follows from
	\eqref{item:convex:reach} and~\cite[4.8\,(12)]{MR0110078}. The second
	equation in \eqref{item:convex:normal} follows from
	\cite[I.2.3]{MR567696}.  Finally, \eqref{item:convex:boundary} is
	implied by \cite[1.3.2]{MR3155183}.
\end{Proof}

\begin{theorem} \label{thm:hausdorff_distance}
	Suppose $X = \mathbf R^n \cap \mathbf B (0,1)$, $F$ is the family of
	nonempty closed subsets of $X$ endowed with the Hausdorff metric, and
	$G = F \cap \{ C \with \textup{$C$ is convex} \}$.

	Then, the following four statements hold.
	\begin{enumerate}
		\item \label{item:hausdorff_distance:compactness} The
		families $F$ and $G$ are compact.
		\item \label{item:hausdorff_distance:dist} The function
		mapping $(x,B) \in X \times F$ onto $\dist (x,B) \in \mathbf
		R$ is continuous.
		\item \label{item:hausdorff_distance:dim} The function
		mapping $C \in G$ onto $\dim C \in \mathbf Z$ is lower
		semicontinuous.
		\item \label{item:hausdorff_distance:boundary} If $\Phi = ( G
		\times F ) \cap \{ (C,B) \with \textup{$B$ is the relative
		boundary of $C$} \}$, then $\Phi$ is a Borel function whose
		domain equals the Borel set $G \cap \{ C \with \dim C \geq 1
		\}$.
	\end{enumerate}
\end{theorem}

\begin{Proof}
	\eqref{item:hausdorff_distance:compactness} is contained in
	\cite[2.10.21]{MR41:1976}.  \eqref{item:hausdorff_distance:dist}
	follows from~\ref{remark:distance}.  We observe that, in order to
	prove
	\eqref{item:hausdorff_distance:dim}~and~\eqref{item:hausdorff_distance:boundary},
	it sufficient to establish the following assertion.  \emph{If $k$ is
	an integer, $C_i$ is a sequence in $G$ with $\dim C_i = k$, $C \in G$,
	and $C_i \to C$ as $i \to \infty$, then $\dim C \leq k$ and, in case
	of equality with $k\geq 1$, also $\Phi(C) = \lim_{i \to \infty} \Phi
	(C_i)$.}  For this purpose, we assume, possibly passing to a
	subsequence, that for some affine subspace~$Q$ of~$\mathbf R^n$
	\begin{equation*}
		\dist (v,\aff C_i) \to \dist (v,Q) \quad \text{as $i \to
		\infty$ for $v \in \mathbf R^n$},
	\end{equation*}
	and, if $k \geq 1$, that for some $B \in F$, we have $\Phi (C_i) \to
	B$ as $i \to \infty$.  It follows $C \subset Q$, whence we
	infer $\dim C \leq \dim Q \leq k$.  Therefore, if $\dim C = k \geq 1$,
	then $Q = \aff C$ and we could assume, possibly replacing $C_i$ by
	$g_i^{-1} [ C_i ]$ for a sequence of isometries~$g_i$ of~$\mathbf R^n$
	with $\lim_{i \to \infty} g_i(x)=x$ for $x \in \mathbf R^n$ and using
	\ref{remark:distance}, that $C_i \subset Q$ for each index~$i$; in
	which case $\Phi (C) = B$ follows readily from
	\ref{lemma:convex}\,\eqref{item:convex:boundary}.
\end{Proof}

\begin{remark} \label{remark:borel}
	We observe that
	\eqref{item:hausdorff_distance:dist}--\eqref{item:hausdorff_distance:boundary}
	imply that, \emph{if $A$ is a Borel subset of $\mathbf R^n$ and
	$\Gamma : A \to G$ is a Borel function, then the set of $(a,v) \in A
	\times \mathbf R^n$ such that $v$~belongs to the relative interior of
	$\Gamma (a)$ is a Borel subset of $\mathbf R^n \times \mathbf R^n$.}
\end{remark}

The corollary to the next theorem on convex cones will be one of the
ingredients to the \hyperlink{observation}{geometric observation for convex
sets} described in the introduction.

\begin{definition} \label{def:cone}
	A subset $C$ of $\mathbf R^n$ is said to be a \emph{cone} if and only
	if $\lambda c \in C$ whenever $0 < \lambda < \infty$ and $c \in C$.
\end{definition}

\begin{theorem} \label{thm:cone_control}
	Suppose $C$ is a convex cone in $\mathbf R^n$,
	\begin{gather*}
		D = \mathbf R^n \cap \{ d \with \textup{$d \bullet c \leq 0$
		for $c \in C$} \},
	\end{gather*}
	$U$ is an $m$~dimensional plane in $\mathbf R^n$, $U \subset D$, $\dim
	C \geq n-m$, and $v$ belongs to the relative interior of $C$.

	Then, $\dim C = n-m$ and there exists $0 \leq \gamma < \infty$
	satisfying
	\begin{equation*}
		\dist (d,U) \leq - \gamma d \bullet v \quad \text{for $d \in
		D$}.
	\end{equation*}
\end{theorem}

\begin{Proof}
	Defining $V = \mathbf R^n \cap \{ v \with \textup{$u \bullet v = 0$
	for $u \in U$} \}$, we see $C \subset V$ from~\cite[4.5]{MR0110078},
	hence $\aff C = V$; in particular, $\dim C = n-m$.  Since $D$ is
	closed under addition and $U \subset D$, $D$ is invariant under
	directions in $U$.  Therefore, it is sufficient to prove the existence
	of $0 \leq \gamma < \infty$ such that the inequality holds for $d \in
	D \cap V \cap \mathbf S^{n-1}$.  If there were no such $\gamma$, then,
	by compactness, there would exist $d \in D \cap V \cap \mathbf
	S^{n-1}$ with $d \bullet v = 0$ which would imply that $v$ belongs to
	the relative boundary of $C$, as $d \in \aff C$.
\end{Proof}

\begin{corollary} \label{corollary:cone_control}
	Under the hypotheses of \ref{thm:cone_control}, there holds
	\begin{equation*}
		\dist (b,U) \leq - \gamma b \bullet v + (1 + \gamma |v|) \dist
		(b,D) \quad \text{for $b \in \mathbf R^n$}.
	\end{equation*}
\end{corollary}

\begin{Proof}
	In view of \ref{remark:distance} and
	\ref{lemma:convex}\,\eqref{item:convex:reach}, one may apply
	\ref{thm:cone_control} to $d = \boldsymbol \xi_D (b)$.
\end{Proof}

\section{Distance bundle}

In the present section, we introduce the distance bundle in
\ref{def:distance_bundle}--\ref{remark:hug}; its nonzero directions correspond
to the \emph{normal bundle} employed by Hug, Last, and Weil
in~\cite{MR2031455}, see~\ref{remark:hug}.  Then, we extend
(in~\ref{thm:properties_dis}) some basic estimates from Federer's treatment of
sets of positive reach in~\cite{MR0110078} which lead to an important
one-sided estimate for the nearest point projection
in~\ref{corollary:one_sided_estimate}.  Finally, we derive the
\hyperlink{observation}{geometric observation}, described for convex sets in
the introduction, in~\ref{lemma:relative_interior}, and the main
\hyperlink{structural}{structural theorem on the singularities of closed sets}
in~\ref{thm:stratification-borel}.

\begin{definition} \label{def:distance_bundle}
	Suppose $A \subset \mathbf R^n$.  Then, the \emph{distance bundle of
	$A$} is defined by
	\begin{equation*}
		\Dis (A) = ( \mathbf R^n \times \mathbf R^n ) \cap \{ (a,v)
		\with \text{$a \in \Clos A$ and $|v| = \dist (a+v,A)$} \}.
	\end{equation*}
	Moreover, we let $\Dis (A,a) = \{ v \with (a,v) \in \Dis (A) \}$ for
	$a \in \mathbf R^n$.
\end{definition}

\begin{remark} \label{remark:dis}
	Clearly, $\Dis (A) = \Dis ( \Clos (A))$, $\Dis (A)$ is closed, and $0
	\in \Dis(A,a)$ if and only if $a \in \Clos A$.  Moreover, $\Dis (A,a)$
	is a convex subset of $\Nor (A,a)$ for $a \in \mathbf R^n$ by
	\cite[4.8\,(2)]{MR0110078}.
\end{remark}

\begin{remark} \label{remark:dis_borel}
	If $X$ and $G$ are as in \ref{thm:hausdorff_distance}, then the
	function mapping $a \in \Clos A$ onto $\Dis(A,a) \cap X \in G$ is a
	Borel function; in fact, \ref{remark:dis}~implies that, in the
	terminology of \cite[II.20]{MR0467310}, the function in question is an
	upper semicontinuous multifunction, whence the assertion follows
	by~\cite[III.3]{MR0467310}.
\end{remark}

\begin{remark} \label{remark:dis_xi}
	If $a \in A$, $v \in \Dis (A,a)$, and $0 \leq t < 1$, then
	$\boldsymbol \xi_A (a+tv)=a$.  In particular, $\boldsymbol \xi_A (a+v)
	= a$ whenever $v$ belongs to the relative interior of $\Dis (A,a)$,
	and $\Dis (A,a)$ is the closure of $\{ v \with \boldsymbol \xi_A (a+v)
	= a \}$.
\end{remark}

\begin{remark} \label{remark:stacho}
	In view of \ref{remark:dis_xi}, we could have alternatively formulated
	our main theorem (see~\ref{thm:stratification-borel}), for closed
	sets, in terms of the bundle $\{ (a,v) \with \boldsymbol \xi_A (a+v) =
	a \}$ which would be more in line with Stach{\'o}'s definition of
	\emph{prenormals} in~\cite[p.~192]{MR534512}.  Our choice of bundle is
	motivated by the fact that $\Dis (A)$ is closed.
\end{remark}

\begin{remark} \label{remark:hug}
	If $A$ is closed, then \ref{remark:dis_xi} yields that
	\begin{equation*}
		\big \{ (a,|v|^{-1}v) \with \textup{$(a,v) \in \mathbf R^n
		\times \mathbf R^n$ and $0 \neq v \in \Dis (A,a)$} \big \}
	\end{equation*}
	equals the \emph{normal bundle} of $A$ defined
	in~\cite[p.~239]{MR2031455}.
\end{remark}

Basic estimates for the distance bundle are collected in the following
theorem.

\begin{theorem} \label{thm:properties_dis}
	Suppose $A \subset \mathbf R^n$.  Then, the following three statements
	hold.
	\begin{enumerate}
		\item \label{item:properties_dis:angle} If $0 < q < \infty$,
		$a \in \Clos A$, $b \in \Clos A$, $v \in \mathbf R^n$, and
		\begin{equation*}
			\text{either} \quad v=0 \quad \text{or} \quad
			q|v|^{-1}v \in \Dis (A,a),
		\end{equation*}
		then $(b-a) \bullet v \leq (2q)^{-1} |b-a|^2 |v|$.
		\item \label{item:properties_dis:projection} If $0 < r < q <
		\infty$, $x \in \mathbf R^n$, $y \in \mathbf R^n$, $a \in A$,
		$b \in A$, and
		\begin{gather*}
			|x-a|=\dist (x,A) \leq r, \quad |y-b| = \dist (y,A)
			\leq r, \\
			\text{either} \quad x = a \quad \text{or} \quad q
			|x-a|^{-1} ( x-a ) \in \Dis (A,a), \\
			\text{either} \quad y = b \quad \text{or} \quad q
			|y-b|^{-1} ( y-b ) \in \Dis (A,b),
		\end{gather*}
		then $\boldsymbol \xi_A (x) = a$, $\boldsymbol \xi_A (y) = b$,
		and
		\begin{equation*}
			| b-a | \leq q(q-r)^{-1} |y-x|.
		\end{equation*}
		\item \label{item:properties_dis:cone} If $0 < q < \infty$, $a
		\in \Clos A$, $b \in \Clos A$, $C$ is a convex cone in
		$\mathbf R^n$,
		\begin{equation*}
			qv \in \Dis (A,a) \quad \text{whenever $v \in C \cap
			\mathbf S^{n-1}$},
		\end{equation*}
		and $D = \mathbf R^n \cap \{ u \with \textup{$u \bullet v \leq
		0$ for $v \in C$} \}$, then
		\begin{equation*}
			\dist (b-a,D) \leq (2q)^{-1} |b-a|^2.
		\end{equation*}
	\end{enumerate}
\end{theorem}

\begin{Proof}
	To prove \eqref{item:properties_dis:angle}, we assume $v \neq
	0$, let $w = |v|^{-1} v$, and compute
	\begin{gather*}
		|a+qw-b| \geq \dist (a+qw,A) = q, \quad |a-b|^2 + 2q w \bullet
		(a-b) + q^2 \geq q^2, \\
		2qw \bullet (b-a) \leq |b-a|^2, \quad v \bullet (b-a) \leq
		(2q)^{-1} |b-a|^2|v|.
	\end{gather*}

	To prove~\eqref{item:properties_dis:projection}, we notice that $a =
	\boldsymbol \xi_A (x)$ and $b = \boldsymbol \xi_A (y)$ by
	\ref{remark:dis_xi} and infer
	\begin{equation*}
		(b-a) \bullet (x-a) \leq |b-a|^2 r/(2q), \quad (a-b) \bullet
		(y-b) \leq |b-a|^2 r/(2q).
	\end{equation*}
	from applying~\eqref{item:properties_dis:angle} twice; once with $v$
	replaced by $x-a$ and once with $a$, $b$, and~$v$ replaced by $b$,
	$a$, and~$y-b$.  Therefore, we obtain
	\begin{align*}
		|b-a| |y-x| & \geq (b-a) \bullet (y-x) \\
		& = (b-a) \bullet ( (b-a) + (a-x) + (y-b) ) \geq |b-a|^2
		(1-r/q),
	\end{align*}
	whence we infer $|x-y| \geq |a-b| (q-r)/q$.

	To prove~\eqref{item:properties_dis:cone}, we suppose $a = 0$.
	Whenever $v \in C$, we notice that
	\begin{equation*}
		v \bullet b \leq (2q)^{-1} |b|^2 |v|
	\end{equation*}
	by~\eqref{item:properties_dis:angle}, and estimate
	\begin{equation*}
		|b-v|^2 = |b|^2 + |v|^2 - 2 b \bullet v \geq |b|^2 + |v|^2 -
		|b|^2 |v|/q \geq |b|^2 - |b|^4 /(4q^2).
	\end{equation*}
	Consequently, $\dist (b,C)^2 \geq |b|^2 - |b|^4/(4q^2)$
	and~\eqref{item:properties_dis:cone} is implied
	by~\cite[4.16]{MR0110078}.
\end{Proof}

\begin{remark}
	The proof is almost verbatim taken from~\cite[4.8\,(7)\,(8),
	4.18\,(2)]{MR0110078}, where slightly stronger hypotheses were made.
\end{remark}

Next, we derive a crucial one-sided estimate for the nearest point projection.

\begin{corollary} \label{corollary:one_sided_estimate}
	Suppose $A \subset \mathbf R^n$, $0 < s < r < q < \infty$, and
	\begin{gather*}
		x \in \dmn \boldsymbol \xi_A, \quad s \leq \dist (x,A) \leq r,
		\quad v = \frac{ x-\boldsymbol \xi_A (x)}{|x - \boldsymbol
		\xi_A (x)|}, \quad q v \in \Dis (A,\boldsymbol \xi_A(x)), \\
		y \in \dmn \boldsymbol \xi_A, \quad s \leq \dist (y,A) \leq r,
		\quad w = \frac{y-\boldsymbol \xi_A (y)}{|y - \boldsymbol
		\xi_A (y)|}, \quad q w \in \Dis (A,\boldsymbol \xi_A(y)).
	\end{gather*}

	Then, there holds
	\begin{equation*}
		( \boldsymbol \xi_A (x) - \boldsymbol \xi_A (y) ) \bullet v
		\leq \kappa |y-x|^2,
	\end{equation*}
	where $\kappa = (2s)^{-1} ( 1 + 2q/(q-r) )^2$.
\end{corollary}

\begin{Proof}
	We let $a = \boldsymbol \xi_A (x)$ and $b = \boldsymbol \xi_A (y)$,
	hence $a = x -|x-a|v$ and $b = y - |y-b|w$.  Next, we
	estimate
	\begin{equation*}
		( a - b ) \bullet v \leq (2s)^{-1} |y-x|^2
	\end{equation*}
	in case $\dist (x,A) = \dist (y,A) = s$; in fact, noting $\dist (y,A)
	\leq |y-a|$ and $|v|=|w|=1$, we obtain
	\begin{gather*}
		s^2 \leq |y-(x-sv)|^2, \quad (x-y) \bullet v \leq (2s)^{-1}
		|y-x|^2, \quad
		(w-v) \bullet v \leq 0, \\
		( a - b ) \bullet v = (x-y) \bullet v + s ( w-v ) \bullet v
		\leq (2s)^{-1} |y-x|^2.
	\end{gather*}

	In the general case, we let
	(see~\ref{lemma:convex}\,\eqref{item:convex:reach})
	\begin{equation*}
		\alpha = a + \boldsymbol \xi_{\mathbf B(0,s)} (x- a), \quad
		\beta = b + \boldsymbol \xi_{\mathbf B(0,s)} ( y-b ),
	\end{equation*}
	notice $\alpha = a +sv$ and $\beta = b + sw$, and infer
	\begin{gather*}
		\alpha \in \dmn \boldsymbol \xi_A, \quad \boldsymbol \xi_A
		(\alpha) = a, \quad \beta \in \dmn \boldsymbol \xi_A, \quad
		\boldsymbol \xi_A (\beta) = b, \\
		|\beta-\alpha| \leq |y-x| + 2 |b-a| \leq ( 1 + 2q/(q-r) )
		|y-x|
	\end{gather*}
	from~\ref{remark:dis_xi},
	\ref{lemma:convex}\,\eqref{item:convex:reach}, and
	\ref{thm:properties_dis}\,\eqref{item:properties_dis:projection}.
	Therefore, we may apply the previous case with $x$ and $y$ replaced by
	$\alpha$ and $\beta$ to deduce the conclusion.
\end{Proof}

\begin{remark}
	One could also derive a two-sided estimate; in fact, this is done in
	the submitted PhD thesis of the second author.
\end{remark}

We now have all ingredients at our disposal to derive the geometric observation, formulated in the introduction
for convex sets, in full generality.

\begin{lemma} \label{lemma:relative_interior}
	Suppose $A \subset \mathbf R^n$, $0 < q < \infty$, $m$ is an integer,
	$1 \leq m < n$, $W$~is the set of $y \in \dmn \boldsymbol \xi_A$
	satisfying
	\begin{equation*}
		0 < \dist (y,A) < q \quad \text{and} \quad  q|y-\boldsymbol
		\xi_A (y)|^{-1} (y-\boldsymbol \xi_A (y)) \in \Dis (A,
		\boldsymbol \xi_A (y)),
	\end{equation*}
	$x \in W$, $a = \boldsymbol \xi_A (x)$, $\dim \Dis (A, a) \geq n-m$,
	$U$ is an $m$~dimensional subspace of~$\mathbf R^n$, $U \subset \Tan (
	A, a)$, and
	\begin{equation*}
		\text{$q|x-a|^{-1}(x-a)$ belongs to the relative interior
		of~$\Dis ( A, a )$}.
	\end{equation*}
	
	Then,
	\begin{equation*}
		\limsup_{W \owns y \to x} |y-x|^{-2} \dist (\boldsymbol \xi_A
		(y)- a,U) < \infty.
	\end{equation*}
\end{lemma}

\begin{Proof}
	Assume $a = 0$, choose $s$ and $r$ such that $0 < s < |x| < r < q$,
	and let $Q = \aff \Dis (A,0)$.  Then, the set~$X$ of all $v \in Q
	\without \{ 0 \}$, such that $q |v|^{-1} v$ belongs to the relative
	interior of $\Dis (A,0)$, is relatively open in~$Q$ and $x \in X$.
	This implies the existence of $\varepsilon > 0$ such that the convex
	cone
	\begin{equation*}
		C = Q \cap \{ v \with \textup{$|rv-x| < \varepsilon$ for some
		$0 < r < \infty$} \}
	\end{equation*}
	satisfies $C \cap \{ v \with |v| = |x| \} \subset X$, hence
	\begin{equation*}
		q v \in \Dis (A,0) \quad \text{whenever $v \in C \cap
		\mathbf S^{n-1}$};
	\end{equation*}
	in particular, $C \subset \Nor (A,0)$ by~\ref{remark:dis}.  We note
	that $\dim C = \dim Q \geq n-m$ and that $x$ belongs to the relative
	interior of $C$, as $Q \cap \mathbf U (x,\varepsilon) \subset C$.
	Abbreviating $D = \mathbf R^n \cap \{ d \with \textup{$d \bullet c
	\leq 0$ for $c \in C$} \}$, we observe $U \subset D$
	from~\cite[4.5]{MR0110078}, and employing $0 \leq \gamma < \infty$
	from~\ref{thm:cone_control} with $v = x$, we estimate
	\begin{align*}
		\dist ( \boldsymbol \xi_A(y),U) & \leq - \gamma
		\boldsymbol \xi_A(y) \bullet x + (1+\gamma|x|) \dist
		(\boldsymbol \xi_A (y), D) \\
		& \leq \gamma \kappa |x| |y-x|^2 + ( 1 + \gamma|x| )
		(2q)^{-1}| \boldsymbol \xi_A (y) |^2 \leq \lambda |y-x|^2
	\end{align*}
	whenever $y \in W$ and $s \leq \dist (y,A) \leq r$ by
	\ref{corollary:cone_control}, \ref{corollary:one_sided_estimate},
	\ref{thm:properties_dis}\,\eqref{item:properties_dis:cone}, and
	\ref{thm:properties_dis}\,\eqref{item:properties_dis:projection},
	where $\kappa = (2s)^{-1} ( 1 + 2 q/(q-r) )^2$ and $\lambda = \gamma
	\kappa |x| + ( 1 + \gamma |x| ) 2^{-1} q (q-r)^{-2}$.  Finally,
	$x$~belongs to the interior of $W \cap \{ y \with s \leq \dist (y,A)
	\leq r \}$ relative to $W$ by \ref{remark:distance}.
\end{Proof}

Finally, we establish the structural theorem on the
singularities of closed sets; in fact, we may formulate it for arbitrary
subsets of Euclidean space.

\begin{theorem} \label{thm:stratification-borel}
	Suppose $A \subset \mathbf R^n$, $m$ is an integer, and $0 \leq m \leq
	n$.  Then,
	\begin{equation*}
		\{ a \with \dim \Dis (A,a) \geq n-m \}
	\end{equation*}
	is a countably $m$ rectifiable Borel set which can be $\mathscr H^m$
	almost covered by the union of a countable family of $m$~dimensional
	submanifolds of~$\mathbf R^n$ of class~$2$.
\end{theorem}

\begin{Proof}
	Let $B = \{ a \with \dim \Dis (A,a) \geq n-m \}$.  We assume $A$ to be
	a nonempty closed set by~\ref{remark:dis}, and also $m < n$. As
	$0 \in \Dis (A,a)$ for $a \in A$ by~\ref{remark:dis}, we obtain
	\begin{equation*}
		\dim \Dis (A,a) = \dim ( \Dis (A,a) \cap \mathbf B
		(0,1)) \quad \text{for $a \in A$}
	\end{equation*}
	from~\ref{lemma:convex}\,\eqref{item:convex:tangent}; in particular,
	$B$~is a Borel set by
	\ref{thm:hausdorff_distance}\,\eqref{item:hausdorff_distance:dim} and
	\ref{remark:dis_borel}.  We define $N$ to be the set of all $(a,v) \in
	A \times \mathbf R^n$ such that $v$ belongs to the relative interior
	of $\Dis (A,a) \cap \mathbf B(0,1)$, hence $N$ is a Borel set by
	\ref{remark:borel} and \ref{remark:dis_borel}.
	By~\ref{remark:dis_xi}, we have
	\begin{equation*}
		\boldsymbol \xi_A (x+v) = a \quad \text{whenever $(a,v) \in
		N$}.
	\end{equation*}
	Noting~\ref{remark:distance} and \ref{remark:xi}, we define~$W_i$ to
	be the Borel set of all $x \in \boldsymbol \xi_A^{-1} [B]$ satisfying
	\begin{equation*}
		0 < \dist (x,A) < i^{-1} \quad \text{and} \quad \big
		(\boldsymbol \xi_A (x), i^{-1} |x-\boldsymbol \xi_A(x)|^{-1} (
		x-\boldsymbol \xi_A (x)) \big ) \in N
	\end{equation*}
	for every positive integer~$i$.  Then,
	$\boldsymbol \xi_A|W_i$ is locally Lipschitzian
	by \ref{thm:properties_dis}\,\eqref{item:properties_dis:projection}
	and \ref{remark:distance},
	and
	\begin{equation*}
		(\boldsymbol \xi_A(x),x-\boldsymbol \xi_A(x)) \in N \quad
		\text{for $x \in W_i$}
	\end{equation*}
	by \ref{remark:relative_interior}.  We observe that this implies that
	\begin{equation*}
		\mathscr H^{n-m} \big ( ( \boldsymbol \xi_A |
		W_i )^{-1} [ \{ \boldsymbol \xi_A (x) \} ] \big ) > 0 \quad
		\text{whenever $x \in W_i$},
	\end{equation*}
	since $( \boldsymbol \xi_A | W_i )^{-1} [ \{ \boldsymbol \xi_A (x) \}
	]$ is relatively open in $\{ \boldsymbol \xi_A (x) + v \with v \in \aff \Dis (A,
	\boldsymbol \xi_A (x)) \}$.

	We choose a countable family $F$ of $m$ dimensional affine planes
	in~$\mathbf R^n$ such that $Q \cap \bigcup F$ is dense in $Q$,
	whenever $Q$ is an affine subspace of $\mathbf R^n$ with $\dim Q \geq
	n-m$; in fact, one may take $F$ to be a countable dense subset in the
	family of all $m$ dimensional affine planes in $\mathbf R^n$.  Thence,
	we deduce, employing \ref{remark:relative_interior}, that
	\begin{equation*}
		{\textstyle B = \bigcup_{i=1}^\infty \boldsymbol \xi_A \left [
		W_i \cap \bigcup F \right ]};
	\end{equation*}
	in fact, whenever $a \in B$, we take $Q = \{ a+v \with v \in \aff
	\Dis(A,a) \}$, pick a positive integer~$i$ such that, for some $x \in
	Q$ with $0 < |x-a| < i^{-1}$, we have that $i^{-1} |x-a|^{-1}
	(x-a)$ belongs to the relative interior of $\Dis (A,a) \cap \mathbf B
	(0,1)$, choose such $x$ within $\bigcup F$, and conclude $x \in W_i$
	with $\boldsymbol \xi_A (x)=a$, as $(a,x-a) \in N$.  It follows that
	$B$ is countably $m$~rectifiable.
	
	To prove the remaining property of $B$, we assume $m \geq 1$.  Then,
	in view of~\ref{remark:coarea_tangent} and
	\ref{lemma:relative_interior}, we may apply
	\ref{corollary:rectifiability_criterion} with $f = \boldsymbol \xi_A |
	W_i$ for every positive integer~$i$ to obtain the conclusion.
\end{Proof}

\begin{remark}
	Our proof of the countable $m$~rectifiability follows
	\cite[4.15\,(3)]{MR0110078}, where the case of sets of positive reach
	was treated.
\end{remark}

\begin{remark} \label{remark:alberti}
	If $A$ is a closed convex set, this property was proven, by
	different methods, in~\cite[Theorem~3]{MR1384392}; the agreement, in
	this case, of the normal bundle used there with our distance bundle
	follows
	from~\ref{lemma:convex}\,\eqref{item:convex:reach}\,\eqref{item:convex:normal}
	and \ref{remark:dis_xi}.
\end{remark}

\begin{remark}
	For $1 \leq m < n$, the preceding theorem may not be strengthened by
	replacing the distance bundle by the normal bundle, as is evident from
	considering a closed $m$~dimensional submanifold of~$\mathbf R^n$ of
	class~$1$ that meets every $m$~dimensional submanifold of~$\mathbf
	R^n$ of class~$2$ in a set of $\mathscr H^m$ measure zero; the
	existence of such $A$ follows from~\cite{MR0427559}.
\end{remark}

\medskip \noindent \textsc{Affiliations}

\medskip \noindent Ulrich Menne \smallskip \\
Institut f{\"u}r Mathematik \\
Mathematisch-naturwissenschaftliche Fakult{\"a}t \\
Universit{\"a}t Z{\"u}rich \\
Winterthurerstrasse 190 \\
\textsc{8057 Z{\"u}rich \\ Switzerland}

\medskip

\noindent Mario Santilli \smallskip \\
Max Planck Institute for Gravitational Physics (Albert
Einstein Institute) \newline Am M{\"u}hlen\-berg 1 \newline 14476
\textsc{Golm} \newline \textsc{Germany} \smallskip \\
University of Potsdam \\
Institute for Mathematics \\
OT Golm \newline Karl-Liebknecht-Stra{\ss}e 24--25 \newline 14476
\textsc{Potsdam}
\newline
\textsc{Germany}

\medskip \noindent \textsc{Email addresses}

\medskip \noindent
\href{mailto:Ulrich.Menne@math.uzh.ch}{Ulrich.Menne@math.uzh.ch} \quad
\href{mailto:Mario.Santilli@aei.mpg.de}{Mario.Santilli@aei.mpg.de}


\begin{thebibliography}{HLW04}

\bibitem[Alb94]{MR1384392}
Giovanni Alberti.
\newblock On the structure of singular sets of convex functions.
\newblock {\em Calc. Var. Partial Differential Equations}, 2(1):17--27, 1994.
\newblock URL: \url{http://dx.doi.org/10.1007/BF01234313}.

\bibitem[Alm86]{MR855173}
F.~Almgren.
\newblock Optimal isoperimetric inequalities.
\newblock {\em Indiana Univ. Math. J.}, 35(3):451--547, 1986.
\newblock URL: \url{http://dx.doi.org/10.1512/iumj.1986.35.35028}.

\bibitem[CV77]{MR0467310}
C.~Castaing and M.~Valadier.
\newblock {\em Convex analysis and measurable multifunctions}.
\newblock Lecture Notes in Mathematics, Vol. 580. Springer-Verlag, Berlin-New
  York, 1977.
\newblock URL: \url{http://dx.doi.org/10.1007/BFb0087685}.

\bibitem[Fed59]{MR0110078}
Herbert Federer.
\newblock Curvature measures.
\newblock {\em Trans. Amer. Math. Soc.}, 93:418--491, 1959.
\newblock URL: \url{https://doi.org/10.1090/S0002-9947-1959-0110078-1}.

\bibitem[Fed69]{MR41:1976}
Herbert Federer.
\newblock {\em Geometric measure theory}.
\newblock Die Grundlehren der ma\-the\-ma\-ti\-schen Wissenschaften, Band 153.
  Springer-Verlag New York Inc., New York, 1969.
\newblock URL: \url{http://dx.doi.org/10.1007/978-3-642-62010-2}.

\bibitem[HLW04]{MR2031455}
Daniel Hug, G{\"u}nter Last, and Wolfgang Weil.
\newblock A local {S}teiner-type formula for general closed sets and
  applications.
\newblock {\em Math. Z.}, 246(1-2):237--272, 2004.
\newblock URL: \url{http://dx.doi.org/10.1007/s00209-003-0597-9}.

\bibitem[Kel75]{MR0370454}
John~L. Kelley.
\newblock {\em General topology}.
\newblock Springer-Verlag, New York, 1975.
\newblock Reprint of the 1955 edition [Van Nostrand, Toronto, Ont.], Graduate
  Texts in Mathematics, No. 27.

\bibitem[Koh77]{MR0427559}
Robert~V. Kohn.
\newblock An example concerning approximate differentiation.
\newblock {\em Indiana Univ. Math. J.}, 26(2):393--397, 1977.
\newblock URL: \url{http://www.iumj.indiana.edu/docs/26030/26030.asp}.

\bibitem[KS80]{MR567696}
David Kinderlehrer and Guido Stampacchia.
\newblock {\em An introduction to variational inequalities and their
  applications}, volume~88 of {\em Pure and Applied Mathematics}.
\newblock Academic Press Inc. [Harcourt Brace Jovanovich Publishers], New York,
  1980.

\bibitem[RZ01]{MR1846894}
J.~Rataj and M.~Z{\"a}hle.
\newblock Curvatures and currents for unions of sets with positive reach. {II}.
\newblock {\em Ann. Global Anal. Geom.}, 20(1):1--21, 2001.
\newblock URL: \url{http://dx.doi.org/10.1023/A:1010624214933}.

\bibitem[San17]{arXiv:1701.07286v1}
Mario Santilli.
\newblock Rectifiability and approximate differentiability of higher order for
  sets, 2017.
\newblock \href {http://arxiv.org/abs/1701.07286v1}
  {\path{arXiv:1701.07286v1}}.

\bibitem[Sch09]{MR2472179}
Reiner Sch{\"a}tzle.
\newblock Lower semicontinuity of the {W}illmore functional for currents.
\newblock {\em J. Differential Geom.}, 81(2):437--456, 2009.
\newblock URL:
  \url{http://projecteuclid.org/getRecord?id=euclid.jdg/1231856266}.

\bibitem[Sch14]{MR3155183}
Rolf Schneider.
\newblock {\em Convex bodies: the {B}runn-{M}inkowski theory}, volume 151 of
  {\em Encyclopedia of Mathematics and its Applications}.
\newblock Cambridge University Press, Cambridge, expanded edition, 2014.
\newblock URL: \url{https://doi.org/10.1017/CBO9781139003858}.

\bibitem[Sta79]{MR534512}
L.~L. Stach{\'o}.
\newblock On curvature measures.
\newblock {\em Acta Sci. Math. (Szeged)}, 41(1-2):191--207, 1979.

\bibitem[Z{\"a}h86]{MR849863}
M.~Z{\"a}hle.
\newblock Integral and current representation of {F}ederer's curvature
  measures.
\newblock {\em Arch. Math. (Basel)}, 46(6):557--567, 1986.
\newblock URL: \url{http://dx.doi.org/10.1007/BF01195026}.

\end{thebibliography}
\end{document}